\documentclass[a4paper,leqno,11pt]{amsart}

\usepackage{amsfonts,amssymb,verbatim,amsmath,amsthm,latexsym,textcomp,amscd}
\usepackage{latexsym,amsfonts,amssymb,epsfig,verbatim}
\usepackage{amsmath,amsthm,amssymb,latexsym,graphics,textcomp}
\usepackage{graphicx}
\usepackage{color}
\usepackage{geometry}
\usepackage{layout}
\usepackage{url}
\usepackage{enumerate}
\usepackage[mathscr]{euscript}
\input xy
\xyoption{all}

\usepackage{hyperref}

\theoremstyle{plain}
\newtheorem{theorem}{Theorem}[section]
\newtheorem{thm}[theorem]{Theorem}

\theoremstyle{definition}

\author{Prabir Das Adhikary, Koustav Banerjee, Manosij Ghosh Dastidar }

\email{dipumsc41@gmail.com}
\address{Pondicherry University,R. V. Nagar, Kalapet, Puducherry, PIN-605014, India.}

\email{banerjeekoustav98@gmail.com; koustav.banerjee@rkmvu.ac.in  }
\address{Ramakrishna Mission Vivekananda University,Belur, Howrah - 711202, West Bengal, India.}

\email{gdmanosij@gmail.com}
\address{Pondicherry University,R. V. Nagar, Kalapet, Puducherry, PIN-605014, India.}

 \keywords {Partitions; cubic partitions; partition congruences.}

\begin{document}

\title[cubic partition $a(n)$]{Some Results On The Cubic Partition}
\maketitle

\begin{abstract}
In this paper we explore Kruyswijk's method and show how to obtain congruences for cubic partition. That apart we also examine inequalities for $a(n)$ and provide upper bound for it in the fashion of the classic partition function $p(n)$.

\end{abstract}

\section{Introduction }

There are many astonishing and beautiful results attributed to Ramanujan.
However, the congruence identities he discovered for the partition
function were the ones described by Hardy as his \lq{most beautiful identity}\rq.
Ramanujan showed that $p(5n+4)\equiv0(mod\thinspace5)$ and this was the first
congruence result regarding the partition function. Over the course
of years many such congruences were discovered. It has also been proven
that such congruence results exists \cite{Ahl-3} for all integers co-prime to 6.
Like the classic partition function $p(n)$ the function $a(n)$ also
known as the cubic partition gives the number of partitions of n where
every even number comes with two different labels and is given by
a similar generating function viz.\begin{eqnarray*}
\sum_{n=0}^{\infty}a(n)x^{n}=\prod_{i=1}^{\infty}\frac{1}{(1-x^{i}){(1-x^{2i})}}.
\end{eqnarray*}
In his paper  \cite{Chan-1}, Chan showed that $a(3n+2)\equiv0(mod\thinspace3)$ which was the
first of such congruence results for the function $a(n)$. In section 3
 we provide an alternative proof of this result. That
apart we exhibit a simple analogous inequality for $a(n)$ as well as
upper bound for it. Here we prove the following theorems:

\begin{thm}
Alternative proof of, $a(3n+2)\equiv0(mod\thinspace3)$.
\end{thm}

\begin{thm}\label{Thm}
 Given any positive integer $n\geq2$,$$a(n+2)+a(n-2)>2a(n)$$
\end{thm}

\begin{thm}
for any positive integer $n$,$$a(n)<\exp^{k\sqrt{n}}$$
where $k=\pi$.
\end{thm}

In the next section, we present some background material on partition function $p(n)$ and some useful theorems due to Euler which will help us to prove our theorems.

\section{Preliminaries }

A partition of a natural number $n$ is a finite sequence of non-increasing positive integer parts where the sum of the parts is equal to $n$. If $p(n)$ denote the number of partitions of $n$, then
$$\sum_{n=0}^\infty p(n)q^n=\dfrac{1}{(q;q)_\infty },$$
where, as customary, for any complex number $a$ and $|q|<1$
$$(a;q)_\infty =\prod_{k=0}^\infty (1-aq^k).$$

Next define integers $a(n)$ by
\begin{equation*}
\sum_{n=0}^\infty a(n)q^n=\dfrac{1}{(q;q)_\infty (q^2;q^2)_\infty }\\
\end{equation*}
Note that $a(n)$ is the number of cubic partitions of $n$ where one of the colors appears only in parts that are multiples of $2$.\\
Here the celebrated Jacobi's triple product identity which states that:\\
For $z\neq 0$ and $\mid q\mid<1$,
$$\sum_{n=-\infty}^\infty {z^{n}q^{n^{2}}={(-zq;q^2)_\infty}{(-q/z;q^2)_\infty}{(q^2;q^2)_\infty}}$$
Now we introduce the Ramanujan's general theta function which is defined by $$f(a,b):= \sum_{n=-\infty}^\infty {a^{n(n+1)/2} b^{n(n-1)/2}} , \mid ab\mid<1$$\\
Three special cases are defined by, in Ramanujan's notation and using Jacobi's Triple product identity we have\\
\begin{equation} \label{eqn I}
\phi(q):= f(q,q) = {(-q,q^2)^2_\infty}{(q^2,q^2)_\infty}
\end{equation}
\begin{equation} \label{eqn II}
\psi(q):= f(q,q^3) = \frac{{(q^2,q^2)_\infty}}{{(q,q^2)^2_\infty}}
\end{equation}

We will be using the following two well-known results due to Euler to prove Theorem 1.1:

\subsection{Result 1 (Generating function for integer partition)}

For $|x|<1$ we have,
$$\prod_{i=1}^\infty \frac{1}{1-x^i} =\sum_{i=0}^\infty p(n)x^{n}$$
where $p(0)=1$.

\vspace{.5cm}
Another very well referred result in integer partition is Euler's pentagonal-number theorem,
which states the following.

\subsection{Result 2 (Euler's pentagonal-number theorem)}\label{result 2}

If $|x|<1$  we  have,

\begin{eqnarray*}
\prod_{m=1}^\infty (1-x^{m}) & = &  1-x-x^{2}+x^{5}+x^{7}-x^{12}-x^{15}+................\\
 & = &  1+\sum_{n=1}^\infty (-1)^{n}({x^{\omega(n)}+x^{\omega(-n)}}) \\
 & = &  \sum_{n=-\infty}^\infty {x^{\omega(n)}}\\ 
\end{eqnarray*}
 where $\omega(n)=\frac {3n^{2}+n}{2}$

\section{Proof of Theorems}
\subsection{Proof of Theorem 1.1}
Chan showed the beautiful identitiy regarding cubic partition $a(n)$\thinspace i.e.\thinspace $a(3n+2)\equiv0(mod\thinspace3)$ \thinspace using modular forms and Roger-Ramanujan continued fraction. In 1950, D.Kruyswijk \cite{Kru-2} gave an alternative proof of celebrated partition congruence result $p(5n+4)\equiv0(mod\thinspace5)$. One can apply the method of Kruyswijk to a certain extent to provide an alternative proof of  $a(3n+2)\equiv0(mod\thinspace3)$.
Here we show explicitly the generating function of $a(3n+2)$ i.e. we establish the following identity:
\begin{eqnarray*}
\sum_{m=0}^{\infty}a(3m+2)x^{m}=3\frac{(E(x^{3}))^{3}(E(x^{6}))^{3}}{E(x)^{4}E(x^{2})^{4}}\\
\end{eqnarray*}
 where $x$ is a complex number satisfying $|x|<1$; 
from the above identity it is clear that $a(3n+2)\equiv0(mod3)$. Where $$E(x)= \prod_{n=1}^{\infty} (1-x^{n})$$\\
\vspace{.5cm}\\
In order to prove the Theorem 1.1 we will prove following lemmas:
\subsubsection{Lemma}
\begin{eqnarray*}
\prod_{h=1}^k(1-x\epsilon^{h})(1-x^{2}\epsilon^{h})& = & (1-x^k)(1-x^{2k})\\
\end{eqnarray*}
where $\epsilon=\exp^{\frac{2\pi i}{k}}$, $i=\sqrt{-1}$
\vspace{.5cm}\\
Proof: To prove the lemma it suffices to prove the identity 
\begin{eqnarray*}
\prod_{h=1}^k(1-x\epsilon^{h})& = & (1-x^k)\\
\end{eqnarray*}
 and one can deduce the lemma from the fact.
 \\ We know that
 \begin{eqnarray*}
 (x^k-1)=\prod_{h=1}^{k}(x-\epsilon^{h})\\
 \end{eqnarray*}
   and so

\begin{eqnarray*}
(1-x^{k}) & = & (-1)^{k+1}\prod_{h=1}^k(\epsilon^{h}-x)\\ & = & (-1)^{k+1}\prod_{h=1}^k(\epsilon^{-h}-x)\\ & = & (-1)^{k+1}\prod_{h=1}^k\epsilon^{-h}\prod_{h=1}^k(1-x\epsilon^{h})\\ & = & (-1)^{k+1}\epsilon^{-\frac{k(k+1)}{2}}\prod_{h=1}^k(1-x\epsilon^{h})\\ & = & \prod_{h=1}^k(1-x\epsilon^{h}).\\
\end{eqnarray*}

 \subsubsection{Lemma} For $gcd(n,k)=d$;

\begin{eqnarray*}
\prod_{h=1}^k(1-x\epsilon^{nh})(1-x^{2}\epsilon^{nh}) & = & (1-x^{\frac{k}{d}})^{d}(1-x^{\frac{2k}{d}})^{d}\\
\end{eqnarray*}
and can deduce that

\begin{eqnarray*}
\prod_{h=1}^k(1-x^{n}\epsilon^{nh})(1-x^{2n}\epsilon^{nh}) & = &  (1-x^{nk})(1-x^{2nk}),\enspace if \thinspace \thinspace gcd(n,k)=1\\
 & = & (1-x^{n})^{k}(1-x^{2n})^{k},\enspace if \thinspace \thinspace k\mid n\\
\end{eqnarray*}
 
Proof: We just compute for $$\prod_{h=1}^k (1-x\epsilon^{nh})$$ and can comment about the product.\\
let $d=gcd(n,k)$ ; $m = \frac{n}{d}$ and let $\delta=\exp^{\frac{2\pi id}{k}}$, $i=\sqrt{-1}$\\
then\enspace
\begin{eqnarray*}
\prod_{h=1}^k(1-x\epsilon^{nh})=\prod_{h=1}^k(1-x\delta^{mh})^{d}
\end{eqnarray*}

since $(m,\frac{k}{d})=1$, $mh$ runs through a complete system of residues mod k,\\
hence by what we have done before\\

\begin{eqnarray*}
\prod_{h=1}^\frac{k}{d}(1-x\delta^{mn})^{d} & = & \prod_{h=1}^\frac{k}{d}(1-x\delta^{h})^{d}\\
 & = & (1-x^{\frac{k}{d}})^{d}
\end{eqnarray*}
so, if $gcd(n,k)=1$, then\\

\begin{eqnarray*}
\prod_{h=1}^{k}(1-x^{n}\exp^{\frac{2\pi inh}{k}}) & = & (1-x^{nk})\\
\end{eqnarray*}

 Also if $k \mid n$ then $gcd(n,k)=k$, so\\

\begin{eqnarray*} 
\prod_{h=1}^{k}(1-x^{n}\exp^{\frac{2\pi inh}{k}}) & = & (1-x^{n})^{k}\\
 \end{eqnarray*}
 And similarly one can do for \\

\begin{eqnarray*}
\prod_{h=1}^k(1-x^{2n}\epsilon^{nh}) & = &  (1-x^{2nk}),\enspace if \thinspace \thinspace gcd(n,k)=1\\
 & = & (1-x^{2n})^{k},\enspace if \thinspace \thinspace k\mid n\\
\end{eqnarray*}
so, we are done.

 \subsubsection{Lemma} 
for prime $q$ and $|x|<1$ we have
\begin{eqnarray*} 
\prod_{n=1}^\infty \prod_{h=1}^q(1-x^{n}\exp^{\frac{2\pi inh}{q}})(1-x^{2n}\exp^{\frac{2\pi inh}{q}}) & = & \frac{E(x^{q})^{q+1}}{E(x^{q^{2}})}\frac{E(x^{2q})^{q+1}}{E(x^{2q^{2}})}\\
\end{eqnarray*}
\vspace{.5cm}\\
where $$E(x)= \prod_{n=1}^{\infty} (1-x^{n})$$\\

Proof:

\begin{eqnarray*}
& &\prod_{n=1}^\infty \prod_{h=1}^q(1-x^{n}\exp^{\frac{2\pi inh}{q}})\prod_{n=1}^\infty \prod_{h=1}^q(1-x^{2n}\exp^{\frac{2\pi inh}{k}})\\
=& &\prod_ {n=1}^\infty(1-x^{qn})^{q}\prod_{r=1}^{q-1}\prod_{m=1}^\infty(1-x^{q(mq-r)})\prod_{n=1}^\infty \prod_{h=1}^q(1-x^{2n}\exp^{\frac{2\pi inh}{q}})\\
=& &E(x^{q})^{q}\prod_{r=0}^{q-1}\prod_{m=1}^\infty (1-x^{q(mq-r)})\prod_{m=1}^\infty(1-x^{q^{2}m})^{-1}\prod_{n=1}^\infty \prod_{h=1}^q(1-x^{2n}\exp^{\frac{2\pi inh}{q}})\\
=& &\frac{E(x^{q})^{q+1}E(x^{2q})^{q+1}}{E(x^{q^{2}})E(x^{2q^{2}})}\\
\end{eqnarray*}

\vspace{.5cm}
Next, we want to give an identity regarding the generating function of $a(n)$ with the help of Lemma 3.

 \subsubsection{Lemma} 
$$ \sum_{m=0}^{\infty}a(m)x^{m}=\frac{E(x^{9})E(x^{18})}{E(x^{3})^{4}E(x^{6})^{4}}
\prod_{n=1}^\infty \prod_{h=1}^2(1-x^{n}\exp^{\frac{2\pi i n h}{3}})(1-x^{2n}\exp^{\frac{2\pi i n h}{3}}) $$\\
\vspace{0.25cm}
Proof: Just apply prime $q=3$ in Lemma 3 we have,\\
\begin{eqnarray*}
\prod_{n=1}^\infty \prod_{h=1}^3(1-x^{n}\exp^{\frac{2\pi i n h}{3}})
(1-x^{2n}\exp^{\frac{2\pi i n h}{3}})=\frac{E(x^{3})^{4}
E(x^{6})^{4}}{E(x^{9})E(x^{18})}\\
\end{eqnarray*}

\begin{eqnarray*}
& &\Rightarrow \prod_{n=1}^\infty \prod_{h=1}^2(1-x^{n}\exp^{\frac{2\pi i n h}{3}})
(1-x^{2n}\exp^{\frac{2\pi inh}{3}})\prod_{n=1}^{\infty}(1-x^{n})(1-x^{2n})=\frac{E(x^{3})^{4}
E(x^{6})^{4}}{E(x^{9})E(x^{18})} \\
\end{eqnarray*}

\begin{eqnarray*}
& &\Rightarrow \frac{E(x^{9})E(x^{18})}{E(x^{3})^{4}E(x^{6})^{4}}
 \prod_{n=1}^\infty \prod_{h=1}^2(1-x^{n}\exp^{\frac{2\pi i n h}{3}})
(1-x^{2n}\exp^{\frac{2\pi i n h}{3}})=\prod_{n=1}^{\infty}\frac{1}{(1-x^{n})(1-x^{2n})}\\
\end{eqnarray*}

\begin{equation} \label{eqn}
\sum_{m=0}^{\infty}a(m)x^{m}=\frac{E(x^{9})E(x^{18})}{E(x^{3})^{4}E(x^{6})^{4}}
\prod_{n=1}^\infty \prod_{h=1}^2(1-x^{n}\exp^{\frac{2\pi i n h}{3}})(1-x^{2n}\exp^{\frac{2\pi inh}{3}})
\end{equation}

\vspace{0.5cm}

Definition: If $q$ is a prime and if $0\leq r<q$, a power series of the form $$\sum_{n=0}^{\infty}b(n)x^{qn+r}$$ is said to be of type $r$ mod $q$.\\
Observe that if $S_{k}$ is a series of type $k$ mod $q$, $S_{m}$ is a series of type $m$ mod $q$ then $S_{k}\times S_{m}$ is a series of type $k+m$ mod $q$.
\vspace{0.5cm}\\
Next, can estimate the product series $$\prod_{n=1}^\infty(1-x^{n})(1-x^{2n})$$ as a sum of power series of type $k mod\thinspace 3$ by Euler pentagonal numbar theorem.\\

 \subsubsection{Lemma}
 We show that
 \begin{eqnarray*} 
 E(x)E(x^{2}) & = & \prod_{n=1}^\infty(1-x^{n})(1-x^{2n})\\ & = & I_{0}+I_{1}+I_{2}\\
 \end{eqnarray*} where $I_{k}$ denote a power series of type\ $k$ mod $3$. 
 \vspace{0.5cm}\\
Proof:\begin{eqnarray*}
& &E(x)E(x^{2})\\
=& &\prod_{n=1}^\infty(1-x^{n})\prod_{n=1}^\infty(1-x^{2n})\\
=& &(1-x-x^{2}+x^{5}+x^{7}-x^{12}-x^{15}+.......)(1-x^{2}-x^{4}+x^{10}+x^{14}-x^{24}-x^{30}+.....)\thinspace by\thinspace \ref{result 2}\\
=& &1-x-2x^{2}+x^{3}+2x^{5}+x^{6}-2x^{9}+x^{10}-2x^{11}-..........\\
=& &(1+x^{3}+x^{6}-2x^{9}-......)+(-x+x^{10}+......)+(-2x^{2}+2x^{5}-......)\\
=& &I_{0}+I_{1}+I_{2}\\
\end{eqnarray*}

In very next lemma we want to estimate the product $$\prod_{n=1}^\infty \prod_{h=1}^2(1-x^{n}\exp^{\frac{2\pi i n h}{3}})(1-x^{2n}\exp^{\frac{2\pi inh}{3}})$$\\
\subsubsection{Lemma}
We show $$\prod_{h=1}^2\prod_{n=1}^\infty(1-x^{n}\alpha^{nh})(1-x^{2n}\alpha^{nh})=\prod_{h=1}^2(I_{0}+I_{1}
\alpha^{h}+I_{2}\alpha^{2h}),\thinspace \thinspace \alpha= \exp^{\frac{2\pi i }{3}}.$$
\vspace{0.5cm}
Proof: Here,
\begin{eqnarray*}
 E(x\alpha^{h})=\prod_{n=1}^\infty(1-x^{n}\alpha^{nh})\\
 \end{eqnarray*}
and

$$ E(x^{2}\alpha^{2h})= \prod_{n=1}^\infty(1-x^{2n}\alpha^{2nh})= 
\prod_{n=1}^\infty(1-x^{2n}\alpha^{nh}) = E(x^{2}\alpha^{h})\\$$

\vspace{0.5cm}
So, now a bit of calculation using Euler pentagonal number theorem we have,
\vspace{0.25cm}\\
\begin{eqnarray*}
E(x\alpha^{h})E(x^{2}\alpha^{h}) & = &\prod_{n=1}^\infty(1-x^{n}\alpha^{nh})(1-x^{2n}\alpha^{nh})\\
& = & (1-x\alpha^{h}-x^{2}\alpha^{2h}+x^{5}\alpha^{2h}+x^{7}\alpha^{h}-x^{12}-x^{15}+...)
(1-x^{2}\alpha^{2h}-x^{4}\alpha^{h}+x^{10}\alpha^{h}+..)\\
& = & (1+x^{3}+x^{6}-2x^{9}-.......)+(-x+x^{10}+.....)\alpha^{h}+(-2x^{2}+2x^{5}-2x^{11}+.....)\alpha^{2h}\\ & = & I_{0}+I_{1}\alpha^{h}+I_{2}\alpha^{2h}\\
\end{eqnarray*}
Therefore,
\begin{eqnarray*}
 \prod_{h=1}^2\prod_{n=1}^\infty(1-x^{n}\alpha^{nh})(1-x^{2n}\alpha^{nh}) & = & \prod_{h=1}^2(I_{0}+I_{1}
\alpha^{h}+I_{2}\alpha^{2h})\\.
\end{eqnarray*}

\vspace{0.5cm}
Then    \ref{eqn}  gives:
\begin{equation}\label{eqn1}
\sum_{m=0}^{\infty}a(m)x^{m}=\frac{E(x^{9})E(x^{18})}{E(x^{3})^{4}E(x^{6})^{4}}\prod_{h=1}^2(I_{0}+I_{1}
\alpha^{h}+I_{2}\alpha^{2h})\\
\end{equation}
So by equating terms of type $2$ mod $3$ in L.H.S and R.H.S of  \ref{eqn1}  we have,
\begin{equation}\label{eqn2}
\sum_{m=0}^{\infty}a(3m+2)x^{3m+2}=B_{2}\frac{E(x^{9})E(x^{18})}{E(x^{3})^{4}E(x^{6})^{4}}\\
\end{equation}\\
Where $B_{2}$ is the power series of type $2$ mod $3$ in $$\prod_{h=1}^2(I_{0}+I_{1}
\alpha^{h}+I_{2}\alpha^{2h})$$\\
\vspace{0.3cm}
( As if $S_{k}$ is a series of type $k$ mod $3$, $S_{m}$ is a series of type $m$ mod $3$ then $S_{k}\times S_{m}$ is a series of type $k+m$ mod $3$ .
 hence $\frac{E(x^{9})E(x^{18})}{E(x^{3})^{4}E(x^{6})^{4}}$ is of type $0$ mod $3$ )\\
 
Next we compute the product $$\prod_{h=1}^2(I_{0}+I_{1}\alpha^{h}+I_{2}\alpha^{2h})$$ to find an explicit series for $B_{2}$.\\
It is clear that 
\begin{eqnarray*}
& &\prod_{h=1}^2(I_{0}+I_{1}\alpha^{h}+I_{2}\alpha^{2h})\\
=& &(I_{0}^{2}-I_{1}I_{2})+(I_{2}^{2}-I_{0}I_{1})+(I_{1}^{2}-I_{2}I_{0})\\
\end{eqnarray*}\\
In above expression $I_{1}^{2}-I_{2}I_{0}$ is the only series of type $2$ mod $3$.\\
So,
\begin{equation}\label{eqn3}
B_{2}=(I_{1}^{2}-I_{2}I_{0}).\\
\end{equation}

\subsubsection{Lemma}
We show $I_{2}I_{0}=-2I_{1}^{2}$ and $I_{1}=-xE(x^{9})E(x^{18})$\\
\vspace{0.5cm}\\

Proof: First we recall the product formula for Ramanujan's general theta function from Ramanujan's lost notebook which is\\
$$f(a,b)f(c,d)= f(ac,bd)f(a^{2}dq,b^{2}cq)+af(acq,bd/q)f(a^{2}dq^{3},b^{2}c/q)+bf(ac/q,bdq)f(a^{2}d/q,b^{2}cq^{3})$$
$E(x)E(x^{2})= I_{0}+I_{1}+I_{2}$\\
Now,\ref{eqn I} and \ref{eqn II} give,\\ 
\begin{eqnarray*}
E(x)E(x^{2}) & = & f(-x,-x^{2})f(-x^{2},-x^{4})\\ & = & f(x^{3},x^{6})f(-x^{9},-x^{9})-xf(x^{3},x^{6})f(-x^{15},-x^{3})-x^{2}f(1,x^{9})f(-x^{15},-x^{3})\\ & = & I_{0}+I_{1}+I_{2}\\
\end{eqnarray*}
Therefore we have,\\
$I_{0} = f(x^{3},x^{6})f(-x^{9},-x^{9})= f(x^{3},x^{6})\phi(-x^{9})$\\
$I_{1} = -xf(x^{3},x^{6})f(-x^{15},-x^{3})$\\
$I_{2} = -x^{2}f(1,x^{9})f(-x^{15},-x^{3})= -2x^{2}\psi(x^{9})f(-x^{15},-x^{3})$\\
Now,
\begin{eqnarray*}
I_{0}I_{2} & = & -2x^{2}f(x^{3},x^{6})f(-x^{15},-x^{3})\phi(-x^{9})\psi(x^{9})\\ & = & -2x^{2}f(x^{3},x^{6})f(-x^{15},-x^{3}){(x^{9};x^{18})_\infty}{(x^{18};x^{18})^{2}_\infty}\\ & = & -2x^{2}f(x^{3},x^{6})f(-x^{15},-x^{3}){(x^{9};x^{9})_\infty}{(x^{18};x^{18})^{2}_\infty}\\
\end{eqnarray*}
Using Jacobi's triple product identity,
$$f(x,y)={(-x;xy)_\infty}{(-y;xy)_\infty}{(xy;xy)_\infty}$$
We have,
\begin{eqnarray*}
f(x^{3},x^{6})f(-x^{3},-x^{15}) & = & {(-x^{3};x^{9})_\infty}{(-x^{6};x^{9})_\infty}{(x^{9};x^{9})_\infty}{(x^{3};x^{18})_\infty}{(x^{15};x^{18})_\infty}{(x^{18};x^{18})_\infty}\\ & = & \frac{{(-x^{3};x^{9})_\infty}{(-x^{6};x^{9})_\infty}{(-x^{9};x^{9})_\infty}}{(-x^{9};x^{9})_\infty} {(x^{9};x^{9})_\infty}{(x^{18};x^{18})_\infty}{(x^{3};x^{18})_\infty}{(x^{15};x^{18})_\infty}\\ & = & {(-x^{3};x^{3})_\infty}{(x^{3};x^{6})_\infty}{(x^{9};x^{9})_\infty}{(x^{18};x^{18})_\infty}\\ & = & {(x^{9};x^{9})_\infty}{(x^{18};x^{18})_\infty}\\
\end{eqnarray*}

Therefore, $I_{1} = -xf(x^{3},x^{6})f(-x^{15},-x^{3}) = -xE(x^{9})E(x^{18})$\\
and so, it is clear that $I_{2}I_{0}=-2I_{1}^{2}$\\

Then  \ref{eqn3} gives, $B_{2}=(I_{1}^{2}-I_{2}I_{0})=3I_{1}^{2} \thinspace$ and $I_{1}=-xE(x^{9})E(x^{18})$\\
 Therefore from  \ref{eqn2} we have \\
 \begin{eqnarray*}
 \sum_{m=0}^{\infty}a(3m+2)x^{3m+2} & = & 3x^{2}\frac{(E(x^{9}))^{3}(E(x^{18}))^{3}}{E(x^{3})^{4}E(x^{6})^{4}}\\
 \end{eqnarray*}
 \vspace{0.5cm}
 Next we dividing both side by $x^{2}$ and replacing $x^{3}$ by $x$, We have\\
 \begin{eqnarray*}
 \sum_{m=0}^{\infty}a(3m+2)x^{m} & = & 3\frac{(E(x^{3}))^{3}(E(x^{6}))^{3}}{E(x)^{4}E(x^{2})^{4}}\\
 \end{eqnarray*}
 Hence $$a(3n+2)\equiv0(mod\thinspace3)$$\\
 
 \subsection{Proof of Theorem 1.2} Here we want to establish a recursive type of inequality of cubic partition of $n$ i.e. $a(n)$ .\\ 
Let $a(n)$ denote the  cubic partitions of $n$ ,i.e. in which even parts of $n$ come up with two colours. Here we denote the two colours as labeling $1$ and $2$ ,viz 2 comes up with $2_{1}$ and $2_{2}$.\\
  We have to proof the given condition  \ref{Thm}  \thinspace $a(n+2)+a(n-2)>2a(n)$.\\
  This condition may be arranged to give $a(n+2)-a(n)>a(n)-a(n-2)$.\\
  Consequently, let's have a look at the cubic partitions of integers $n$\thinspace and\thinspace$n+2$. Now, if $2_{1}$ is added as an exrtra part at the end of a cubic partition of n ,a cubic partition of $n+2$ is obtained. In particular a cubic partition that contains $2_{1}$.\\
  Conversely, if $2_{1}$ is deleted from a cubic partition of $n+2$, which contains a $2_{1}$,a cubic partition of n results. Thus there is one to one correspondence between the entire set of cubic partition of n, and the subset of cubic partitions of $n+2$ that contain  $2_{1}$. Hence $a(n+2)-a(n)$ is the number of cubic partitions of  $n+2$ in the subset X of cubic partitions which do not contain $2_{1}$. Similarly $a(n)-a(n-2)$ is the number of cubic partition of $n$ in the subset Y of cubic partitions which do not contain $2_{1}$.\\
 If we can show that cardinality of X greater than cardinality  of Y, then the desired conclusion follows. This simply requires we show that, to each overpartition in Y , there corresponds a distinct overpartition in X.\\
 Now we define a one-one function from Y to X, which is not onto.\\
 let $f:Y\longrightarrow X$ be defined by\\
  $y_{1}+y_{2}+.....+y_{k}\longrightarrow y_{1}+y_{2}+...+(2+y_{i})+..+y_{k} $ , $y_{i}=max(y_{1},y_{2},.....y_{k})$, where an even number (not equal to 2)with two different label we prefer the even number with label 1,(viz $4_{1}>4_{2}$)\\
  Clearly this is well defined one-one function but not onto (since 1+1+....+1 is a member of X whose pre image is not in Y).\\
 Thus, cardinality of X greater than cardinality of Y and we are done.\\

\subsection{Proof of Theorem 1.3} In this section we give an upper bound for $a(n)$ almost like the classical upper bound for $p(n)$ by estimating the logarithm of generating function $F(x)$ for $a(n)$\\
where  $$F(x)=\prod_{i=1}^{\infty}\frac{1}{(1-x^{i})(1-x^{2i})}$$\\
Now,\\
\begin{eqnarray*}
F(x) & = & \prod_{i=1}^{\infty}\frac{1}{(1-x^{i})(1-x^{2i})}\\ & = & 1+\sum_{m=1}^{\infty}a(n)x^{n}\\
\end{eqnarray*}
and restrict to the interval $0<x<1$.\\
Then clearly $a(n)x^{n}<F(x)$ , from which we obtain, $\log a(n) < \log F(x)+n \log \frac{1}{x}$.\\
We estimate the terms $\log F(x)$ and $n \log \frac{1}{x}$ separately.
First, we note that\\
\begin{eqnarray*}
\log F(x) & = & -\sum_{i=1}^{\infty}\log (1-x^{i})-\sum_{i=1}^{\infty}\log (1-x^{2i})\\
& = & \sum_{m=1}^{\infty}\frac{1}{m}\frac{x^{m}}{1-x^m}+\sum_{m=1}^{\infty}\frac{1}{m}\frac{x^{2m}}{1-x^{2m}}\\
\end{eqnarray*}

Since we have the following identity\\
$$\frac{1-x^{m}}{1-x}=1+x+x^{2}+.....+x^{m-1}$$\\
\vspace{0.3cm}
 for $0<x<1$ one can deduce that\\
 $$mx^{m-1}<\frac{1-x^m}{1-x}<m $$\\
so, $$\frac{1}{m^{2}}\frac{x^{m}}{1-x}\leq\frac{1}{m}\frac{x^{m}}{1-x^{m}}\leq\frac{1}{m^{2}}\frac{x}{1-x}$$\\
and sum over $m$ we have, 
\begin{equation}\label{eqn4}
\sum_{m=1}^{\infty}\frac{1}{m^{2}}\frac{x^{m}}{1-x}\leq\sum_{m=1}^{\infty}\frac{1}{m^{2}}\frac{x}{1-x}=\frac{\pi^{2}}{6t},\thinspace\thinspace where \thinspace\thinspace t =\frac{1-x}{x} \\
\end{equation}

\vspace{0.3cm}
Similarly, for $0<x<1$ one can deduce that\\
\begin{equation}\label{eqn5}
\sum_{m=1}^{\infty}\frac{1}{2m^{2}}\frac{x^{2m}}{1-x}\leq\sum_{m=1}^{\infty}\frac{1}{m}\frac{x^{2m}}{1-x^{2m}}\leq\sum_{m=1}^{\infty}\frac{1}{2m^{2}}\frac{x}{1-x}\\
\end{equation}

\vspace{0.3cm}
From   \ref{eqn4}   and   \ref{eqn5}   we get,\\
$$\log F(x)\leq \frac{x}{1-x}\sum_{m=1}^{\infty}\frac{1}{m^{2}}+\frac{x}{1-x}\sum_{m=1}^{\infty}\frac{1}{2m^{2}}$$\\
We defined $t=\frac{1-x}{x}$. Also note that t varies from $\infty$ to $0$ through positive values as $x$ varies from 0 to 1.\\
\vspace{0.3cm}
Finally, estimate the term \thinspace\thinspace$ n \log \frac{1}{x}$ and give an upper bound for $a(n)$.\\ 
Now, for $t>0$ we have $\log (1+t)<t$.\\
\vspace{0.5cm}
Since $t=\frac{1-x}{x}$, therefore  $\log \frac{1}{x}<t $.\\
\vspace{0.5cm}
Thus, $$\log a(n)< \log F(x)+n \log \frac{1}{x}<\frac{\pi^{2}}{4t}+nt $$\\
In order to get that bound we need to check for which value of $t$ minimum of $\frac{\pi^{2}}{4t}+nt$ occurs.\\
One can easily check that for $t=\frac{\pi}{2\sqrt{n}}$ minimum occurs and for this value of $t$, we have\\
$$\log a(n)<\pi\sqrt{n}$$.\\
Hence,we are done.\\

\section{Remark}
Chen and lin proved another congruence result on cubic partition $a(25n+22)\equiv0(mod\thinspace5)$ using modular forms in their paper  \cite{Che-4}. One may provide our method to prove the above result.

\mbox{ }\\

\end{document}